\documentclass[11pt,reqno]{amsart}
\usepackage{amsmath,amsthm,amssymb}
\usepackage{mathrsfs}
\usepackage{amsfonts}
\usepackage[noadjust]{cite}
\usepackage[textwidth=150mm,textheight=220mm,
centering]{geometry}
\allowdisplaybreaks
\usepackage{xcolor}
\usepackage[colorlinks=true, pdfstartview=FitV, linkcolor=blue, citecolor=blue, urlcolor=blue]{hyperref}

\newcommand{\bbC}{\mathbb{C}}
\newcommand{\bbR}{\mathbb{R}}
\newcommand{\bbT}{\mathbb{T}}
\newcommand{\bbZ}{\mathbb{Z}}
\newcommand{\ee}[1]{\mathrm{e}^{#1}}
\newcommand{\phaseip}[2]{\langle #1,#2\rangle}
\newcommand{\dist}{\operatorname{dist}}
\newcommand{\Wron}{\mathcal{W}_{\pi}}
\newcommand{\rmd}{\mathrm d}
\makeatletter
\newcommand{\resultalias}[2]{%
  \phantomsection
  \begingroup
  \edef\@currentlabel{\getrefnumber{#2}}%
  \label{#1}%
  \endgroup
}
\makeatother

\newtheorem{theorem}{Theorem}[section]
\newtheorem{definition}[theorem]{Definition}
\newtheorem{lemma}[theorem]{Lemma}
\newtheorem{corollary}[theorem]{Corollary}

\newtheorem{proposition}[theorem]{Proposition}
\newtheorem{remark}[theorem]{Remark}

\numberwithin{equation}{section}

\title{Sampling Density for Multi-window STFT Phase Retrieval}

\author[W. Sun]{Wenchang Sun}
\address{Wenchang Sun:
	School of Mathematical Sciences and LPMC\\
	Nankai University\\
	Tianjin, 300071\\
	China}
\email{sunwch@nankai.edu.cn}

\author[M. Zhou]{Mengxuan Zhou}
\address{Mengxuan Zhou:
	School of Mathematical Sciences and LPMC\\
	Nankai University\\
	Tianjin, 300071\\
	China}
\email{zhoumx@mail.nankai.edu.cn}

\keywords{Phase retrieval; short-time Fourier transform; Hermite windows; Bargmann--Fock space; Beurling density}

\thanks{Corresponding author: Wenchang Sun. E-mail: sunwch@nankai.edu.cn}
\date{}

\thanks{This work was partially supported by the
National Natural Science Foundation of China (12571104).}

\subjclass[2020]{94A20}

\begin{document}

\begin{abstract}
We study phase retrieval for multi-window short-time Fourier transforms with windows given by complex linear combinations of the first two normalized Hermite functions.
Assume that the family of window coefficient vectors does phase retrieval in $\bbC^2$. We prove that the multi-window magnitude samples uniquely determine every signal in $L^2(\bbR)$ up to a global phase whenever the sampling set is separated with lower Beurling density greater than $2$, or is a translated lattice with density at least $2$. Moreover, we show that the sampling density is the best possible,
which improves a result in
[P.~Grohs, L.~Liehr, and M.~Rathmair.
Multi-window {STFT} phase retrieval: lattice uniqueness,
J. Funct. Anal., 288(3): Paper No. 110733, 23, 2025].
\end{abstract}

\maketitle

\section{Introduction and Main Results}
\label{sec:intro}


Phase retrieval is concerned with recovering a signal from its intensity measurements,
which is very useful in coherent diffraction imaging
\cite{Barmherzig2021,
EdidinSuresh2026,
LiuChen2026,
Khare2023},
signal processing
\cite{AmirBendory2026,
Balan2006,
Bandeira2014,
HuangWen2026,
Jaming2026,
PohlYang2014,
Qaddura2026,
XiaXu2026}
and quantum mechanics.
And there are many works on this topic,
which include phase retrieval
for bandlimited signals \cite{AW21,Thakur2011,Wellershoff2023},
for functions in shift invariant spaces
\cite{ChenCheng2020,ChengJiang2019,Groechenig2020},
for wavelet transforms,
and for short-time Fourier transforms (STFT).
In this paper, we focus on STFT phase
retrieval.

The short-time Fourier transform   describes the time-frequency structure of a signal through localized Fourier analysis and is a basic tool in time-frequency analysis.

Given a window function $g\in L^2(\bbR)$,
the STFT of a signal  $f$ is defined by
\begin{equation}\label{eq:stft}
V_gf(x,\omega)=\int_{\bbR}f(t)\overline{g(t-x)}\ee{-2\pi it\omega}\rmd t,
\qquad (x,\omega)\in\bbR^2.
\end{equation}
The aim of STFT phase retrieval is to recover $f$ from the magnitudes $|V_gf|$.
Since $V_g(\tau f)=\tau V_gf$ for every
$\tau $
in
\[
   \bbT:=\{z\in\bbC:|z|=1\},
\]
two signals that differ by a global phase have the same magnitude of STFT.
Throughout this paper, uniqueness is therefore understood with respect to the following equivalence relation:
\[
f\sim h
\iff
\exists \tau\in\bbT:=\{z\in\bbC:|z|=1\},f=\tau h
\]

Discrete magnitude sampling differs fundamentally from linear sampling that retains the complex phase.
Grohs and Liehr~\cite{Grohs2022} proved that, for any single window and any lattice, the short-time Fourier magnitude samples on that lattice cannot do phase retrieval on all of $L^2(\bbR)$. And increasing the lattice density alone cannot remove this obstruction.
A natural alternative is to use several windows on the same sampling set.
This approach retains the structure of regular sampling while recovering local phase relations through interference information between different windows.

Grohs, Liehr, and Rathmair~\cite{GrohsLiehrRathmair2025} established a connection between finite-dimensional phase retrieval, multi-window STFT measurements, and uniqueness in Fock spaces.
Before stating their results, we introduce some notations.

A set $\Lambda \subset\bbR^2$ is called separated if the distances between distinct points have a common positive lower bound.
For a locally finite set $\Lambda $, define its lower and upper Beurling densities by
\begin{align*}
D^-(\Lambda )&=\liminf_{r\to\infty}\inf_{(x,y)\in\bbR^2}
   \frac{\#(\Lambda \cap [x-r,x+r]\times [y-r,y+r] )}{4r^2},\\
D^+(\Lambda )&=\limsup_{r\to\infty}\sup_{(x,y)\in\bbR^2}
   \frac{\#(\Lambda \cap [x-r,x+r]\times [y-r,y+r] )}{4r^2}.
\end{align*}
If these densities coincide, their common value is denoted by $D(\Lambda )$.
Translation and complex conjugation preserve these densities.
For $A\in\mathrm{GL}_2(\bbR)$ and $v\in\bbR^2$, the translated lattice $v+A\bbZ^2$ has the density $1/|\det A|$.

Let
\begin{equation}\label{eq:windows}
 h_0(t)=2^{1/4}\ee{-\pi t^2},\quad
h_1(t)=2^{5/4}\sqrt\pi t \ee{-\pi t^2},\quad
\mbox{ and }\quad
\mathcal H_{01} =\operatorname{span}_{\bbC}\{h_0,h_1\}.
\end{equation}
For $p=(p_1,p_2)^T\in\bbC^2$, denote $g_p=p_1h_0+p_2h_1$.

Given a family of coefficient vectors $\mathcal P\subset\bbC^2$ and a sampling set $\Lambda\subset\bbR^2$,
we say that  $(\{g_p\}_{p\in\mathcal P},\Lambda)$ does STFT phase retrieval on  $L^2(\bbR)$
if, for every $f,h\in L^2(\bbR)$,
\[	\left(|V_{g_p}f(\lambda)|=|V_{g_p}h(\lambda)|,
	\quad\forall p\in\mathcal P,\ \forall\lambda\in\Lambda\right)
	\implies f\sim h.
\]

The problem
of multi-window STFT phase retrieval
reduces to studying the injectivity of the measurement map
\[
\mathcal M_{\mathcal P,\Lambda}(f)
:=\left(|V_{g_p}f(\lambda)|\right)_{p\in\mathcal P,\,\lambda\in\Lambda}
\]
on $L^2(\bbR)/{\sim}$.
When $\mathcal P$ does phase retrieval in $\bbC^2$, that is,
for all $u,v\in\bbC^2$,
\[
(|\phaseip{u}{p}|=|\phaseip{v}{p}|,\forall p\in\mathcal P)
\implies u\sim v,
\]
a fundamental problem is to determine
the density of the sampling sequence $\Lambda$ which ensures this injectivity.

Grohs, Liehr, and Rathmair~\cite{GrohsLiehrRathmair2025}
proved  that if
$\mathcal P\subset\bbC^2$ does phase retrieval in $\bbC^2$
and
$A\in\mathrm{GL}_2(\bbR)$
satisfies
\[
   |\det A|\le \frac{1}{4},
\]
then
$(\{g_p\}_{p\in\mathcal P},A\bbZ^2)$ does
STFT phase retrieval on  $L^2(\bbR)$.

In this paper, we show that $|\det A|\le 1/2$ suffices
for
$(\{g_p\}_{p\in\mathcal P},A\bbZ^2)$ to do
STFT phase retrieval on  $L^2(\bbR)$.
And the condition
$|\det A|\le 1/2$ is also necessary
whenever $A = \begin{pmatrix}
\alpha & 0 \\
0 &\beta
\end{pmatrix}$.

Moreover, for a general separate sampling sequence $\Lambda$, we also obtain the optimal
sampling density. Specifically,
we have the following result.

\begin{theorem}\label{thm:lattice}
Suppose that  $\mathcal P\subset\bbC^2$ does phase retrieval in $\bbC^2$.
Let
$\Lambda=v+A\bbZ^2$ for
some $v\in\bbR^2$ and  $A\in\mathrm{GL}_2(\bbR)$.
If $ |\det A|\le 1/2$, then
$(\{g_p\}_{p\in\mathcal P},\Lambda)$ does phase retrieval on  $L^2(\bbR)$.

The result is sharp in the sense that
$(\{g_p\}_{p\in\mathcal P},\Lambda)$  does
 not do phase retrieval on  $L^2(\bbR)$
whenever  $A = \begin{pmatrix}
\alpha & 0 \\
0 &\beta
\end{pmatrix}$ with
$|\alpha\beta|> 1/2$.
\end{theorem}

For general separated sampling sets, we obtain the following sufficient density condition.

\begin{theorem}\label{thm:intro-density}
Suppose that  $\mathcal P\subset\bbC^2$ does phase retrieval in $\bbC^2$.
If $\Lambda$ is separated and $D^-(\Lambda)>2$,
then $(\{g_p\}_{p\in\mathcal P},\Lambda)$ does phase retrieval on  $L^2(\bbR)$.
\end{theorem}

To prove the main results, we investigate the sampling theory for Fock spaces, which provides the analytic basis for these density questions \cite{Groechenig2001,Zhu2012,Seip1992}.
Uniqueness for lattices at the critical density goes back to Perelomov~\cite{Perelomov1971}.
For the connection between Hermite windows and polyanalytic Fock spaces, we refer to \cite{Abreu2012}.
The modified $\sigma$-function used here comes from the lattice construction of Gr\"ochenig and Lyubarskii~\cite{Groechenig2009}.
We use products associated with finitely many lattice translates to establish uniqueness.
For nonuniqueness, we combine the modified $\sigma$-function with elliptic functions to construct counterexamples on rectangular lattices.

The paper is organized as follows. 
Section~\ref{sec:prelim} introduces the notation and collects the required facts from time-frequency analysis.
Section~\ref{sec:uniqueness} establishes a criterion linking local measurements to uniqueness for Wronskians.
Section~\ref{sec:periodic} studies periodic zero sets and derives explicit density conditions.
Section~\ref{sec:sharp} establishes optimality on rectangular lattices.

\section{Preliminaries and notations}\label{sec:prelim}

In this section, we collect some preliminaries
and notations to be used in the proof.

For convenience, we identify $(x,\omega)\in\bbR^2$ with $x+i\omega\in\bbC$.

Moyal's identity~\cite[Corollary 3.2.2]{Groechenig2001} is
useful in the study of  STFT.

\begin{proposition}\label{lem:basic}
	For $f,g\in L^2(\bbR)$,
	\begin{equation}\label{eq:moyal}
		\|V_gf\|_{2}=\|f\|_2\|g\|_2.
	\end{equation}
\end{proposition}

For $S\subset\bbC$, write $\overline S=\{\overline z:z\in S\}$. This notation denotes the complex conjugate of the set, not its closure.

We write $\rmd\mu$ for planar Lebesgue measure, $B_r(w)=\{z:|z-w|<r\}$, and $\dist(z,S)=\inf_{w\in S}|z-w|$.

\begin{definition}
Let $\mathcal A$ be a class of functions defined on $\Omega$.
A set $Z\subset\Omega$ is called a uniqueness set for $\mathcal A$ if $F\in\mathcal A$ and $F|_Z=0$ imply $F\equiv0$.
\end{definition}

\subsection{Bargmann–Fock spaces}
For $\alpha>0$ and $0<p<\infty$, define
the Bargmann–Fock space as
\[
  \mathcal F^p_\alpha=\{F\text{ entire}:\|F\|_{p,\alpha}<\infty\},
\]
where
\[
\|F\|_{p,\alpha}^p
=\frac{\alpha p}{2\pi}\int_{\bbC}|F(z)|^p \ee{-p\alpha|z|^2/2}\rmd\mu(z).
\]

For $F\in\mathcal F^p_\alpha$,
we have the following pointwise estimate \cite[Corollary 2.8]{Zhu2012},

\begin{equation}\label{eq:pointwise}
|F(z)|\le\|F\|_{p,\alpha}\ee{\alpha|z|^2/2},\qquad z\in\bbC.
\end{equation}

For $f\in L^2(\bbR)$, its Bargmann transform is defined by
\begin{equation}\label{eq:bargmann}
\mathcal Bf(z)=2^{1/4}\int_{\bbR}f(t)\ee{2\pi tz-\pi t^2-\pi z^2/2}\rmd t.
\end{equation}
The Bargmann transform is a unitary isomorphism from $L^2(\bbR)$ onto $\mathcal F^2_\pi$~\cite[Theorem 3.4.3]{Groechenig2001}.

\subsection{Hermite functions}
The normalized Hermite functions can be written as
\[
h_n(t)=\frac{(-1)^n2^{1/4-n}}{\sqrt{\pi^n n!}}
\ee{\pi t^2}\frac{\mathrm d^n}{\mathrm dt^n}\ee{-2\pi t^2},\qquad n\ge 0,
\]
which satisfy $\|h_n\|_2=1$.

For $z=x+iy$, set $\eta(z)=\ee{-i\pi xy+\pi|z|^2/2}$. We see from \cite[Equation~(2.6)]{Abreu2012}
that
\[
\eta(z)V_{h_n}f(\overline z)
=\frac{1}{\sqrt{\pi^n n!}}
\sum_{j=0}^n\binom nj(-\pi\overline z)^{n-j}(\mathcal Bf)^{(j)}(z).
\]
In particular, if $F=\mathcal Bf$, then
\begin{equation}\label{eq:hermite}
F(z)=\eta(z)V_{h_0}f(\overline z)
\mbox{\quad and \quad}
F'(z)=\eta(z)\left(\pi\overline z V_{h_0}f(\overline z)+\sqrt\pi V_{h_1}f(\overline z)\right).
\end{equation}
The complex conjugation in \eqref{eq:hermite} comes from the sign convention in the Fourier exponential in \eqref{eq:stft}.
Thus, the entire-function sampling set corresponding to a time-frequency sampling set $\Lambda$ is $\overline\Lambda$.

\subsection{Stable sampling sets}
A separated set $\Gamma$ is called a stable sampling set for $\mathcal F^2_\alpha$ if there exist $c,C>0$ such that
\[
c\|F\|_{2,\alpha}^2\le
\sum_{z\in \Gamma}|F(z)|^2\ee{-\alpha|z|^2}
\le C\|F\|_{2,\alpha}^2,
\qquad \forall F\in\mathcal F^2_\alpha.
\]
Every stable sampling set is a uniqueness set.

\begin{proposition}\label{lem:sampling}
Let $\Gamma\subset\bbC$ be separated. Then $\Gamma$ is a stable sampling set for $\mathcal F^2_\alpha$ if and only if
\[
D^-(\Gamma)>\frac{\alpha}{\pi}.
\]
\end{proposition}
This classical result can be found in \cite[Theorem~2.1]{Seip1992} and \cite[Corollary~4.37]{Zhu2012}.
The strict inequality is necessary for stable sampling. To treat the critical endpoint for lattices and finite unions of lattice translates, we will prove uniqueness separately rather than infer it from stable sampling.

\section{Wronskian estimates and a phase retrieval criterion}\label{sec:uniqueness}

In this section, we first recover a two-dimensional vector at each sampling point, and then use uniqueness for entire functions to replace the point-dependent phases by a single global phase.
The key step is to control the growth of the Wronskian.

\begin{lemma}\label{lem:wronskian}
	For arbitrary $f,h\in L^2(\bbR)$, write
	\[
	F=\mathcal{B}f,\qquad
	H=\mathcal{B}h,\qquad
	W=FH'-F'H.
	\]
	Then
	\[
	W\in\mathcal{F}^1_{2\pi}\cap\mathcal{F}^2_{2\pi},\qquad
	\|W\|_{1,2\pi}\le 2\sqrt{\pi}\|f\|_2\|h\|_2,
	\qquad
	\|W\|_{2,2\pi}\le 2\sqrt{2\pi}\|f\|_2\|h\|_2.
	\]
\end{lemma}

\begin{proof}
	By \eqref{eq:hermite},
	\[
	F(z) = \eta(z) V_{h_0}f(\overline{z}),\qquad
	F'(z) = \eta(z)(\pi\overline{z}  V_{h_0}f(\overline{z}) + \sqrt{\pi}V_{h_1}f(\overline{z})),
	\]
	\[
	H(z) = \eta(z) V_{h_0}h(\overline{z}),\qquad
	H'(z) = \eta(z)(\pi\overline{z} V_{h_0}h(\overline{z}) + \sqrt{\pi}V_{h_1}h(\overline{z})).
	\]
	Substituting these expressions into $W$ gives
	\[
	W(z)
	= \sqrt{\pi}(\eta(z))^2(V_{h_0}f(\overline{z})V_{h_1}h(\overline{z})-V_{h_0}h(\overline{z})V_{h_1}f(\overline{z})).
	\]
	Since $|\eta(z)|^2=\ee{\pi|z|^2}$, we obtain
	\[
	\|W\|_{1,2\pi}
	=\sqrt{\pi}\|V_{h_0}fV_{h_1}h-V_{h_0}hV_{h_1}f\|_1.
	\]
	Minkowski's inequality gives
	\[
	\|W\|_{1,2\pi}
	\le\sqrt{\pi}\|V_{h_0}fV_{h_1}h\|_1
	+\sqrt{\pi}\|V_{h_0}hV_{h_1}f\|_1.
	\]
	By H\"older's inequality and Proposition~\ref{lem:basic}, together with $\|h_0\|_2=\|h_1\|_2=1$, we have
	\begin{align*}
	\|W\|_{1,2\pi}
	&\le\sqrt{\pi}\|V_{h_0}f\|_2\|V_{h_1}h\|_2
	+\sqrt{\pi}\|V_{h_0}h\|_2\|V_{h_1}f\|_2\\
	&=2\sqrt{\pi}\|f\|_2\|h\|_2.
	\end{align*}
	Similarly, the definition of the Fock norm yields
	\[
	\|W\|_{2,2\pi}
	=\sqrt{2\pi}\|V_{h_0}fV_{h_1}h-V_{h_0}hV_{h_1}f\|_2.
	\]
	Applying Minkowski's inequality and then H\"older's inequality gives
	\begin{align*}
	\|W\|_{2,2\pi}
	&\le\sqrt{2\pi}\|V_{h_0}fV_{h_1}h\|_2
	+\sqrt{2\pi}\|V_{h_0}hV_{h_1}f\|_2\\
	&\le\sqrt{2\pi}\|V_{h_0}f\|_\infty\|V_{h_1}h\|_2
	+\sqrt{2\pi}\|V_{h_0}h\|_\infty\|V_{h_1}f\|_2.
	\end{align*}
	By \eqref{eq:hermite} and \eqref{eq:pointwise},
	\[
	\|V_{h_0}f\|_\infty\le\|F\|_{2,\pi}=\|f\|_2,
	\qquad
	\|V_{h_0}h\|_\infty\le\|H\|_{2,\pi}=\|h\|_2.
	\]
	Proposition~\ref{lem:basic} therefore gives
	\[
	\|W\|_{2,2\pi}\le2\sqrt{2\pi}\|f\|_2\|h\|_2.
	\]
	Thus $W\in\mathcal{F}^1_{2\pi}\cap\mathcal{F}^2_{2\pi}$.
\end{proof}

Write
\[
\Wron=\{FH'-F'H:F,H\in\mathcal F^2_\pi\}.
\]
Here $\Wron$ is used only as a class of entire functions. It is not required to be a linear space.

\begin{lemma}\label{prop:abstract}
	Let $\mathcal P\subset\bbC^2$ do phase retrieval.
	If $\overline\Lambda$ is a uniqueness set for both $\mathcal F^2_\pi$ and $\Wron$, then
	$(\{g_p\}_{p\in\mathcal P},\Lambda)$ does phase retrieval in $L^2(\bbR)$.
\end{lemma}

\begin{proof}
	Suppose that $f,h\in L^2(\bbR)$ satisfy
	\begin{equation}\label{eq:local-hy}
		|V_{g_p}f(z)|=|V_{g_p}h(z)|,\qquad \forall p\in\mathcal{P},\forall z\in\Lambda.
	\end{equation}
	Write
	\[
	F=\mathcal{B}f,\qquad
	H=\mathcal{B}h,\qquad
	W=FH'-F'H.
	\]
	For every $z\in\Lambda$, $f\in L^2(\bbR)$, and $p=(p_1,p_2)\in\bbC^2$, we have
	\[
	|V_{g_p}f(z)|
	= |\overline{p_1} V_{h_0}f(z)+\overline{p_2} V_{h_1}f(z)|
	= \left|\left\langle\begin{pmatrix}V_{h_0}f(z)\\V_{h_1}f(z)
\end{pmatrix},p\right\rangle\right|.
	\]
	Since $\mathcal{P}$ does phase retrieval, \eqref{eq:local-hy} implies
	\[
	\begin{pmatrix}
		V_{h_0}f(z) \\
		V_{h_1}f(z)
	\end{pmatrix}\sim
	\begin{pmatrix}
		V_{h_0}h(z) \\
		V_{h_1}h(z)
	\end{pmatrix},\qquad
	\forall z\in\Lambda.
	\]
	Together with \eqref{eq:hermite}, this yields
	\[
	\begin{pmatrix}F(z)\\F'(z)\end{pmatrix}\sim\begin{pmatrix}H(z)\\H'(z)\end{pmatrix},\qquad\forall z\in\overline\Lambda.
	\]
	Moreover, $W$ vanishes on $\overline{\Lambda}$. The assumption on $\Lambda$ implies that $W$ vanishes identically on $\bbC$.
	
	If $F$ vanishes identically on $\bbC$, then $H$ vanishes on $\overline\Lambda$. The assumption on $\overline\Lambda$ implies that $H$ vanishes identically on $\bbC$. The injectivity of $\mathcal{B}$ gives $f=h=0$, and in particular $f\sim h$.
	
	If $F$ does not vanish identically on $\bbC$, there is a domain $E\subset\bbC$ on which $F$ has no zeros. For every $z\in E$,
	\[
	(H/F)'(z)=\frac{(FH'-F'H)(z)}{F^2(z)}=0.
	\]
	Hence there exists a constant $c\in\bbC$ such that $H=cF$ on $E$. By the identity theorem for entire functions, $H=cF$ on $\bbC$.
	Moreover, since
	\[
	|F(z)|=|H(z)|=|cF(z)|,\qquad
	\forall z\in\overline\Lambda,
	\]
	the function $(1-|c|)F$ vanishes on $\overline\Lambda$. If $F$ vanished on $\overline\Lambda$, the assumption on $\overline\Lambda$ would imply that $F$ vanished identically on $\bbC$, a contradiction. Thus $F$ cannot vanish on all of $\overline\Lambda$, and consequently $|c|=1$. Since the Bargmann transform is injective, the identity
	\[
	\mathcal{B}(h-cf)=H-cF=0
	\]
	implies $f\sim h$. Therefore $(\{g_p\}_{p\in\mathcal{P}},\Lambda)$ does phase retrieval in $L^2(\bbR)$.
\end{proof}

\begin{theorem}\label{thm:fockcriterion}
Let $\mathcal P\subset\bbC^2$ do phase retrieval.
If $\overline\Lambda$ is a uniqueness set for $\mathcal F^2_{2\pi}$, then
$(\{g_p\}_{p\in\mathcal P},\Lambda)$ does phase retrieval in $L^2(\bbR)$.
\end{theorem}

\begin{proof}
Lemma~\ref{lem:wronskian} gives $\Wron\subset\mathcal F^2_{2\pi}$.
Also, for $F\in\mathcal F^2_\pi$,
\[
\|F\|_{2,2\pi}^2=2\int_{\bbC}|F(z)|^2\ee{-2\pi|z|^2}\rmd\mu(z)
\le2\|F\|_{2,\pi}^2,
\]
so $\mathcal F^2_\pi\subset\mathcal F^2_{2\pi}$.
The conclusion follows from Lemma~\ref{prop:abstract}.
\end{proof}

We can now use this theorem to give a proof of Theorem~\ref{thm:intro-density}.

\begin{proof}[Proof of Theorem~\ref{thm:intro-density}]
	Complex conjugation preserves separation and Beurling densities.
	Thus $\overline\Lambda$ is separated and $D^-(\overline\Lambda)>2$.
	By Proposition~\ref{lem:sampling}, $\overline\Lambda$ is a stable sampling set for $\mathcal F^2_{2\pi}$ and hence a uniqueness set for that space.
	Theorem~\ref{thm:fockcriterion} gives the conclusion.
\end{proof}

Theorem~\ref{thm:fockcriterion} separates the choice of windows from the choice of sampling locations: the former is a phase retrieval problem in $\bbC^2$, while the latter is a uniqueness problem in a Fock space with parameter $2\pi$.
In the next section, we turn this criterion into explicit density conditions and treat the critical endpoint.

\section{Uniqueness and density conditions for periodic sampling sets}\label{sec:periodic}

For general separated sets, Proposition~\ref{lem:sampling} already covers densities strictly above the critical value.
To include the critical endpoint for periodic sampling sets, we construct entire functions with prescribed periodic zero sets.
This construction also prepares for the counterexamples to uniqueness in the next section.

Let $\Lambda=\omega_1\bbZ+\omega_2\bbZ$, where $\Im(\omega_2/\omega_1)>0$.
Write
\[
\Pi_\Lambda=\{t\omega_1+u\omega_2:0\le t,u\le1\},\qquad
s(\Lambda)=\Im(\overline{\omega_1}\omega_2)>0.
\]
Here $\Pi_\Lambda$ is a closed fundamental parallelogram
and $s(\Lambda)$ is its area.
If $\Lambda=A\bbZ^2$, then $s(\Lambda)=|\det A|$.
The Weierstrass $\sigma$-function is defined by
\[
\sigma_\Lambda(z)=z\prod_{\omega\in \Lambda\setminus\{0\}}
\left(1-\frac z\omega\right)
\exp\left(\frac z\omega+\frac{z^2}{2\omega^2}\right).
\]
This product converges locally uniformly on compact sets. $\sigma_\Lambda$ is an odd entire function whose zeros are exactly $\Lambda$, and all of them are simple.
There exist quasiperiod constants $\eta_{\Lambda,1},\eta_{\Lambda,2}$ such that
\[
\sigma_\Lambda(z+\omega_k)=-\ee{\eta_{\Lambda,k}(z+\omega_k/2)}\sigma_\Lambda(z),\qquad k=1,2.
\]
The periods $\omega_k$ and constants $\eta_{\Lambda,k}$ also satisfy the Legendre relation
\[
\eta_{\Lambda,1}\omega_2-\eta_{\Lambda,2}\omega_1=2\pi i.
\]
These basic facts about elliptic functions can be found in \cite[Page~150]{Lawden1989} and \cite[Page~15]{Zhu2012}.

The following proposition summarizes the modified Weierstrass
$\sigma$-function construction from
\cite[Propositions 3.4 and 3.5]{Groechenig2009}.

\begin{proposition}\label{lem:sigma}
	Let $\Lambda=\omega_1\bbZ+\omega_2\bbZ$, where $\Im(\omega_2/\omega_1)>0$. Set
	\[
	\kappa_\Lambda = \frac{\eta_{\Lambda,2}\overline{\omega_1}-\eta_{\Lambda,1}\overline{\omega_2}}{\omega_1\overline{\omega_2}-\omega_2\overline{\omega_1}},\qquad
	\gamma_\Lambda = \frac{\pi}{s(\Lambda)},
	\]
	and
	\[
	\nu_\Lambda(z)
	= \sigma_\Lambda(z)\ee{\kappa_\Lambda z^2/2},\qquad
	\rho_\Lambda(z)
	= |\nu_\Lambda(z)|\ee{-\gamma_\Lambda|z|^2/2}.
	\]
	Then
	\begin{enumerate}
		\item $\nu_\Lambda$ is entire and has simple zeros precisely at the points of the lattice $\Lambda$.
		\item $\rho_\Lambda$ is continuous and nonnegative, with zero set exactly $\Lambda$.
		\item $\rho_\Lambda$ is $\Lambda$-periodic, and
		\[
		\rho_\Lambda(z)
		\le \sup_{w\in\Pi_\Lambda}\rho_\Lambda(w),\qquad
		\forall z\in\bbC.
		\]
		\item For every $\varepsilon>0$, there exists a constant $c_\varepsilon>0$ such that $\operatorname{dist}(z,\Lambda)\ge\varepsilon$ implies $\rho_\Lambda(z)\ge c_\varepsilon$.
	\end{enumerate}
\end{proposition}

\begin{theorem}\label{thm:periodic}
	Let $\Lambda=A\mathbb{Z}^2$ be a lattice, where $A\in\mathrm{GL}_2(\mathbb{R})$.
	Suppose that $v_1,\cdots,v_n\in\mathbb{R}^2$ satisfy
	\[
	v_j-v_k\notin\Lambda,
	\qquad 1\le j<k\le n.
	\]
	For each $j=1,\cdots,n$, write
	\[
	\Lambda_j=v_j+\Lambda,
	\qquad
	\Gamma=\bigcup_{j=1}^{n}\Lambda_j.
	\]
	Then $\Gamma$ is a uniqueness set for
	$\mathcal{F}_\alpha^2$ if and only if
	\[
	|\det A|\le\frac{n\pi}{\alpha}.
	\]
\end{theorem}

\begin{proof}
	Write
	\[
	\gamma_\Lambda=\frac{\pi}{|\det A|},
	\qquad
	\gamma=n\gamma_\Lambda
	=\frac{n\pi}{|\det A|}.
	\]
	It suffices to prove that $\Gamma$ is a uniqueness set for $\mathcal{F}_\alpha^2$ if and only if $\gamma\ge\alpha$.
	
	\medskip
	\noindent
	\textbf{Step 1: Construction of an entire function with zero set $\Gamma$.}
\nopagebreak
	
	By Proposition~\ref{lem:sigma}, there exists a nonzero entire function $\nu_\Lambda$ with simple zeros precisely at the points of $\Lambda$, such that
	\[
	\rho_\Lambda(z)
	=|\nu_\Lambda(z)|\ee{-\gamma_\Lambda|z|^2/2}
	\]
	is a continuous, nonnegative, $\Lambda$-periodic function with zero set exactly $\Lambda$. Periodicity and compactness of the fundamental domain imply that there exists a constant $M>0$ such that
	\[
	0\le\rho_\Lambda(z)\le M,
	\qquad \forall z\in\bbC.
	\]
	Moreover, for every $\varepsilon>0$, there exists a constant $c_\varepsilon>0$ such that $\operatorname{dist}(z,\Lambda)\ge\varepsilon$ implies $\rho_\Lambda(z)\ge c_\varepsilon$.
	
	For $j=1,\cdots,n$, define
	\[
	\nu_{\Lambda,j}(z)
	=
	\ee{\gamma_\Lambda\overline{v_j}z-\gamma_\Lambda|v_j|^2/2}
	\nu_\Lambda(z-v_j).
	\]
	Clearly, $\nu_{\Lambda,j}$ is a nonzero entire function with simple zeros precisely at the points of $\Lambda_j$. Define also
	\[
	\rho_{\Lambda,j}(z)
	=
	|\nu_{\Lambda,j}(z)|\ee{-\gamma_\Lambda|z|^2/2}.
	\]
	A direct calculation gives
	\begin{align*}
		\rho_{\Lambda,j}(z)
		&=
		|\nu_\Lambda(z-v_j)|
		\ee{\gamma_\Lambda\operatorname{Re}(\overline{v_j}z)
			-\gamma_\Lambda|v_j|^2/2-\gamma_\Lambda|z|^2/2}\\
		&=
		|\nu_\Lambda(z-v_j)|
		\ee{-\gamma_\Lambda|z-v_j|^2/2}\\
		&=
		\rho_\Lambda(z-v_j).
	\end{align*}
	Thus, each $\rho_{\Lambda,j}$ is a continuous, nonnegative,
	$\Lambda$-periodic function with zero set exactly
	$\Lambda_j$, and
	\[
	0\le\rho_{\Lambda,j}(z)\le M,
	\qquad \forall z\in\bbC.
	\]
	
	Define the entire function
	\[
	G(z)=\prod_{j=1}^{n}\nu_{\Lambda,j}(z)
	\]
	and the function
	\[
	R(z)=|G(z)|\ee{-\gamma|z|^2/2}.
	\]
	Since $\gamma=n\gamma_\Lambda$,
	\begin{equation}
		\label{eq:n-trans-R-product}
		R(z)
		=
		\prod_{j=1}^{n}\rho_{\Lambda,j}(z)
		=
		\prod_{j=1}^{n}\rho_\Lambda(z-v_j).
	\end{equation}
	
	By the assumption $v_j-v_k\notin\Lambda$ ($j\ne k$), the translated lattices $\Lambda_j$ are pairwise disjoint. Thus $G$ has simple zeros precisely at the points of $\Gamma$.
	It follows from \eqref{eq:n-trans-R-product} that
	$R$ is a continuous, nonnegative, $\Lambda$-periodic function
	with zero set exactly $\Gamma$, and
	\begin{equation}
		\label{eq:n-trans-R-upper}
		0\le R(z)\le M^n,
		\qquad \forall z\in\bbC.
	\end{equation}
	For any $\varepsilon>0$,
	if $\operatorname{dist}(z,\Gamma)\ge\varepsilon$,
	then, for all $j=1,\cdots,n$,
	\[
	\operatorname{dist}(z-v_j,\Lambda)
	=
	\operatorname{dist}(z,\Lambda_j)
	\ge\varepsilon.
	\]
	Consequently,
	$\rho_\Lambda(z-v_j)\ge c_\varepsilon$. By \eqref{eq:n-trans-R-product}, whenever $\operatorname{dist}(z,\Gamma)\ge\varepsilon$, we have $R(z)\ge c_\varepsilon^n$.
	Finally, since $R$ is continuous, nonnegative, periodic,
	and not identically zero,
	\begin{equation}
		\int_{\bbC}R(z)^2\rmd\mu(z)
		=\infty.
		\label{eq:n-trans-R-divergence}
	\end{equation}
	
	\medskip
	\noindent
	\textbf{Step 2: If $\gamma<\alpha$, then $\Gamma$ is not a uniqueness set.}
\nopagebreak
	
	The definitions of $G$ and $R$, together with
	\eqref{eq:n-trans-R-upper}, give
	\begin{align*}
		\int_{\bbC}
		|G(z)|^2\ee{-\alpha|z|^2}\rmd\mu(z)
		&=
		\int_{\bbC}
		R(z)^2\ee{-(\alpha-\gamma)|z|^2}\rmd\mu(z)\\
		&\le
		M^{2n}
		\int_{\bbC}
		\ee{-(\alpha-\gamma)|z|^2}\rmd\mu(z)\\
		&=
		\frac{\pi M^{2n}}{\alpha-\gamma}
		<\infty.
	\end{align*}
	Thus $G\in\mathcal{F}_\alpha^2$.
	On the other hand, $G\not\equiv0$
	and $G|_\Gamma=0$.
	Hence there exists a nonzero function
	$G\in\mathcal{F}_\alpha^2$
	that vanishes on $\Gamma$,
	so $\Gamma$ is not a uniqueness set for
	$\mathcal{F}_\alpha^2$.
	
	\medskip
	\noindent
	\textbf{Step 3: If $\gamma\ge\alpha$, then
		$\Gamma$ is a uniqueness set.}
	\nopagebreak

	Let $F\in\mathcal{F}_\alpha^2$
	satisfy $F|_\Gamma=0$.
	We show that $F\equiv0$.
	
	Since the zero set of $G$ is exactly $\Gamma$
	and all its zeros are simple,
	the quotient
	\[
	H(z)=\frac{F(z)}{G(z)},
	\qquad z\in\bbC\setminus\Gamma,
	\]
	has a removable singularity at each point of $\Gamma$.
	Thus $H$ extends uniquely to an entire function
	on $\bbC$.

	For $j\ne k$,
	the assumption $v_j-v_k\notin\Lambda$ and
	the discreteness of $\Lambda$ give
	\[
	\operatorname{dist}(\Lambda_j,\Lambda_k)
	=
	\operatorname{dist}(v_j-v_k,\Lambda)
	>0.
	\]
	Since the number $n$ of translated lattices is finite,
	there exists $\delta>0$
	such that any two distinct points
	$w,w'\in\Gamma$ satisfy $|w-w'|\ge\delta$.
	In other words,
	\[
	\delta_\Gamma
	:=
	\inf_{\substack{w,w'\in\Gamma\\w\ne w'}}
	|w-w'|>0.
	\]
	Fix
	$\varepsilon\in(0,\delta_\Gamma/2)$.
	If $\operatorname{dist}(z,\Gamma)
	\ge\varepsilon$, then
	\eqref{eq:pointwise}
	gives
	\begin{align}
		|H(z)|
		&=
		\frac{|F(z)|}{|G(z)|}
		=
		\frac{|F(z)|}
		{R(z)\ee{\gamma|z|^2/2}}
		\notag\\
		&\le
		\frac{\|F\|_{2,\alpha}}{c_\varepsilon^n}
		\ee{(\alpha-\gamma)|z|^2/2}
		\notag\\
		&\le
		\frac{\|F\|_{2,\alpha}}{c_\varepsilon^n},
		\label{eq:n-trans-h-bound}
	\end{align}
	
	We now extend this estimate to the whole complex plane.
	Since $\varepsilon<\delta_\Gamma/2$,
	the disks $\{B_\varepsilon(w):w\in\Gamma\}$
	are pairwise disjoint.
	Moreover, for every $w\in\Gamma$
	and $z\in\partial B_\varepsilon(w)$,
	we have $\operatorname{dist}(z,\Gamma)
	=\varepsilon$.
	Thus, by \eqref{eq:n-trans-h-bound},
	\[
	|H(z)|
	\le
	\frac{\|F\|_{2,\alpha}}{c_\varepsilon^n},
	\qquad
	\forall z\in\partial B_\varepsilon(w).
	\]
	The maximum modulus principle gives
	\[
	|H(z)|
	\le
	\frac{\|F\|_{2,\alpha}}{c_\varepsilon^n},
	\qquad
	\forall z\in B_\varepsilon(w).
	\]
	Together with \eqref{eq:n-trans-h-bound},
	this yields
	\[
	|H(z)|
	\le
	\frac{\|F\|_{2,\alpha}}{c_\varepsilon^n},
	\qquad
	\forall z\in\bbC.
	\]
	Therefore $H$ is a bounded entire function.
	By Liouville's theorem,
	$H$ is constant.
	Hence there exists $b\in\bbC$
	such that
	\[
	F(z)=bG(z),
	\qquad\forall  z\in\bbC.
	\]
	
	Suppose that $b\ne0$.
	Since $\gamma\ge\alpha$, it follows from
	\eqref{eq:n-trans-R-divergence}
	that
	\begin{align*}
		\int_{\bbC}
		|F(z)|^2\ee{-\alpha|z|^2}\rmd\mu(z)
		&=
		|b|^2
		\int_{\bbC}
		|G(z)|^2\ee{-\alpha|z|^2}\rmd\mu(z)\\
		&=
		|b|^2
		\int_{\bbC}
		R(z)^2\ee{(\gamma-\alpha)|z|^2}\rmd\mu(z)\\
		&\ge
		|b|^2
		\int_{\bbC}R(z)^2\rmd\mu(z)\\
		&=\infty.
	\end{align*}
	This contradicts $F\in\mathcal{F}_\alpha^2$. Therefore $b=0$, and hence $F\equiv0$.
	Since $F$ was arbitrary,
	$\Gamma$ is a uniqueness set for
	$\mathcal{F}_\alpha^2$.
\end{proof}

When $n=1$, Theorem~\ref{thm:periodic} reduces to the classical uniqueness condition for translated lattices.
All translates here must come from a single lattice. The union of two disjoint, distinct lattices need not have a uniform positive separation, so the preceding proof does not apply directly to such a union.

\begin{remark}\label{rem:endpoint}
If $D(\Gamma)=\alpha/\pi$, Theorem~\ref{thm:periodic} guarantees uniqueness, whereas Proposition~\ref{lem:sampling} shows that $\Gamma$ is not a stable sampling set for $\mathcal F^2_\alpha$.
Periodicity allows \eqref{eq:n-trans-R-divergence} to rule out nonzero functions vanishing on the set even at the critical density.
This argument depends on the periodic structure. Density alone does not extend the endpoint conclusion to arbitrary separated sets.
\end{remark}

\begin{proof}[Proof of the sufficient condition in Theorem~\ref{thm:lattice}]
	Complex conjugation maps $\Lambda=v+A\bbZ^2$ to a translated lattice with the same fundamental area $|\det A|\le1/2$.
	Applying Theorem~\ref{thm:periodic} with $n=1$ and Fock parameter $2\pi$ shows that $\overline\Lambda$ is a uniqueness set for $\mathcal F^2_{2\pi}$.
	The conclusion follows from Theorem~\ref{thm:fockcriterion}.
\end{proof}

\begin{corollary}\label{cor:periodic}
	Let $\mathcal P\subset\bbC^2$ do phase retrieval. Let $A\in\mathrm{GL}_2(\bbR)$, and let $v_1,\cdots,v_n\in\bbR^2$ satisfy
	\[
	v_j-v_k\notin A\bbZ^2,\qquad
	1\le j<k\le n.
	\]
	Furthermore, let
	\[
	\Lambda=\bigcup_{j=1}^n(v_j+A\mathbb{Z}^2).
	\]
	If $|\det A|\le n/2$, then
	$(\{g_p\}_{p\in\mathcal P},\Lambda)$ does phase retrieval in $L^2(\bbR)$.
\end{corollary}

\begin{proof}
	The set $\overline\Lambda$ is also a union of $n$ distinct translates of a fixed lattice, with the same density.
	Theorem~\ref{thm:periodic} shows that it is a uniqueness set for $\mathcal F^2_{2\pi}$, and Theorem~\ref{thm:fockcriterion} then gives the result.
\end{proof}

We next give explicit window families satisfying the finite-dimensional condition. The following result is given by \cite[Theorem 1.4]{GrohsLiehrRathmair2025}.

\begin{proposition}\label{thm:four}
	Let $\varphi_0,\varphi_1,\varphi_2,\varphi_3\in\bbC^2$, and define
	\[
	\lambda_k:=\langle\varphi_k,\varphi_0\rangle,\qquad
	\nu_k:=\left\langle\varphi_k,
\left(\begin{smallmatrix}0&-1\\1&0\end{smallmatrix}
\right)\overline{\varphi_0}\right\rangle,\qquad k=1,2,3.
	\]
	The following statements are equivalent:
	\begin{enumerate}
		\item $\{\varphi_0,\varphi_1,\varphi_2,\varphi_3\}\subset\bbC^2$ does phase retrieval.
		\item $\nu_1\nu_2\nu_3\neq0$, and $\lambda_1/\nu_1,\lambda_2/\nu_2,\lambda_3/\nu_3$ are not collinear.
	\end{enumerate}
\end{proposition}

Combining this proposition with the sufficient condition in Theorem~\ref{thm:lattice} and Theorem~\ref{thm:intro-density} immediately yields the following result.

\begin{theorem}
	Let $\lambda_1,\lambda_2,\lambda_3\in\bbC$ be noncollinear, and let
	\[
	\mathcal{P} :=
	\{(1,0)^T,(\lambda_1,1)^T,(\lambda_2,1)^T,(\lambda_3,1)^T\}
	\]
	If $\Lambda\subset\bbR^2$ satisfies one of the following conditions, then
	$(\{g_p\}_{p\in\mathcal P},\Lambda)$ does phase retrieval in $L^2(\bbR)$:
	\begin{enumerate}
		\item $\Lambda$ is separated and $D^-(\Lambda)>2$.
		\item $\Lambda=v+A\bbZ^2$ is a translated lattice and $|\det A|^{-1}\ge2$.
	\end{enumerate}
\end{theorem}

For example, choosing $(\lambda_1,\lambda_2,\lambda_3)=(0,1,i)$ gives
\begin{equation}\label{eq:fourwindows}
h_0,\qquad h_1,\qquad h_0+h_1,\qquad ih_0+h_1.
\end{equation}
The measurements with the last two windows provide, respectively, the real and imaginary parts of the cross term between the first two STFT components.
The four-window construction and the density conditions are logically independent. Our density improvement does not require redesigning this family of windows.

\section{Optimality on rectangular lattices}\label{sec:sharp}

The preceding section established a sufficient condition at density $2$.
To show that this constant cannot be uniformly lowered within the fixed window space $\mathcal H_{01}$, it is not enough to construct a nonzero Fock-space function vanishing on the sampling set: the complete local first-order data of two signals must differ by a common unimodular factor at each sampling point.
Elliptic functions on rectangular lattices provide such a construction.

\begin{proof}[Proof of the necessary condition in Theorem~\ref{thm:lattice}]
	Replacing the spacings by their absolute values leaves the lattice unchanged. We first consider the unshifted rectangular lattice with $\alpha,\beta>0$ and $\alpha\beta>1/2$.
	Write
	\[
	\Lambda=\alpha\bbZ+i\beta\bbZ,\qquad
	\Gamma=2\alpha\bbZ+2i\beta\bbZ,\qquad
	\gamma_\Gamma=\frac{\pi}{s(\Gamma)}
	=\frac{\pi}{4\alpha\beta}<\frac{\pi}{2}.
	\]
	
	Consider the Weierstrass $\wp$-function associated with the lattice $\Gamma$~\cite[Page~13]{Zhu2012},
	\[
	\wp(z)
	=
	\frac{1}{z^2}
	+
	\sum_{\omega\in\Gamma\setminus\{0\}}
	\left(
	\frac{1}{(z-\omega)^2}-\frac{1}{\omega^2}
	\right),\qquad
	z\in\bbC\setminus\Gamma.
	\]
	By \cite[Proposition 1.16]{Zhu2012},
	$\wp$ is an even meromorphic function whose poles are precisely the points of $\Gamma$,
	all of order two. Moreover, $\wp$ is $\Gamma$-periodic.
	Thus $\wp'$ is odd and also $\Gamma$-periodic.
	Since $\overline{\Gamma}=\Gamma$,
	the series representation above gives
	\[
	\wp(\overline z)
	=
	\overline{\wp(z)},\qquad
	\forall z\in\bbC\setminus\Gamma.
	\]
	
	For every $\lambda\in\Lambda\setminus\Gamma$,
	we have $2\lambda\in\Gamma$ and $\overline\lambda-\lambda\in\Gamma$.
	Periodicity implies
	\[
	\overline{\wp(\lambda)}
	=
	\wp(\overline\lambda)
	=
	\wp(\lambda),
	\]
	and
	\[
	\wp'(\lambda)
	=
	\wp'(-\lambda)
	=
	-\wp'(\lambda).
	\]
	Hence
	\[
	\wp(\lambda)\in\bbR,\qquad
	\wp'(\lambda)=0,\qquad
	\forall\lambda\in\Lambda\setminus\Gamma.
	\]
	
	By Proposition~\ref{lem:sigma}, there exists an entire function $\nu_\Gamma$
	with simple zeros precisely at the points of $\Gamma$, such that
	\[
	\rho_\Gamma(z)=|\nu_\Gamma(z)|\ee{-\gamma_\Gamma|z|^2/2}
	\]
	is continuous, nonnegative, and $\Gamma$-periodic.
	Define
	\[
	U(z)=\nu_\Gamma(z)^2,\qquad
	F(z)=U(z)(\wp(z)+i),\qquad
	H(z)=U(z)(\wp(z)-i),
	\qquad z\in\bbC\setminus\Gamma.
	\]
	Since $U$ has a double zero at each $\omega\in\Gamma$
	and $\wp$ has a pole of order two at $\omega$,
	the singularities of $F,H$ at the points of $\Gamma$ are removable.
	We continue to write $F,H$ for their entire extensions to $\bbC$.
	
	We first prove that $F,H\in\mathcal F_\pi^2$.
	For $z\in\bbC\setminus\Gamma$,
	\[
	|F(z)|\ee{-\gamma_\Gamma|z|^2}
	=
	\rho_\Gamma(z)^2|\wp(z)+i|,
	\qquad
	|H(z)|\ee{-\gamma_\Gamma|z|^2}
	=
	\rho_\Gamma(z)^2|\wp(z)-i|.
	\]
	By the periodicity of $\rho_\Gamma$ and $\wp$,
	both right-hand sides are $\Gamma$-periodic.
	Since the left-hand sides extend continuously to $\bbC$,
	their extensions remain $\Gamma$-periodic.
	Periodicity and continuity on a compact fundamental domain imply that there exists a constant $C>0$ such that
	\[
	|F(z)|\le C\ee{\gamma_\Gamma|z|^2},\qquad
	|H(z)|\le C\ee{\gamma_\Gamma|z|^2},\qquad
	\forall z\in\bbC.
	\]
	Since $\pi-2\gamma_\Gamma>0$,
	\[
	\int_\bbC
	\left(|F(z)|^2+|H(z)|^2\right)
	\ee{-\pi|z|^2}\rmd\mu(z)
	\le
	2C^2\int_\bbC\ee{-(\pi-2\gamma_\Gamma)|z|^2}\rmd\mu(z)
	=
	\frac{2\pi C^2}{\pi-2\gamma_\Gamma}
	<\infty.
	\]
	Therefore $F,H\in\mathcal F_\pi^2$.
	
	We next show that $F\not\sim H$.
	If there were $\tau\in\bbT$ such that $F=\tau H$,
	then, since $U$ has no zeros on $\bbC\setminus\Gamma$, we would have
	\[
	\wp(z)+i
	=
	\tau(\wp(z)-i),\qquad
	\forall z\in\bbC\setminus\Gamma.
	\]
	Equivalently,
	\[
	(1-\tau)\wp(z)+i(1+\tau)=0,\qquad
	\forall z\in\bbC\setminus\Gamma.
	\]
	If $\tau=1$, this gives $2i=0$, a contradiction. If $\tau\neq1$, then
	\[
	\wp(z)= -\frac{i(1+\tau)}{1-\tau},\qquad \forall z\in\mathbb C\setminus\Gamma .
	\]
	By the identity theorem for meromorphic functions, this would imply that $\wp$ is constant as a meromorphic function on $\bbC$, which is impossible because $\wp$ has poles of order two at every point of $\Gamma$.
	Thus $F\not\sim H$.
	
	We now prove that, for every $z\in\Lambda$, there exists $\tau_z\in\bbT$ such that
	\[
	F(z) = \tau_z H(z),\qquad
	F'(z) = \tau_z H'(z)
	\]
	First, let $z\in\Lambda\setminus\Gamma$.
	Set
	\[
	r_z=\wp(z)\in\bbR,\qquad
	\tau_z=\frac{r_z+i}{r_z-i}\in\bbT.
	\]
	Since $\wp'(z)=0$,
	\[
	F(z) = (r_z+i) U(z),\quad
	F'(z) = (r_z+i) U'(z),\quad
	H(z) = (r_z-i) U(z),\quad
	H'(z) = (r_z-i) U'(z)
	\]
	Therefore
	\[
	F(z) = \tau_z H(z),\qquad
	F'(z) = \tau_z H'(z)
	\]
	Next, let $z\in\Gamma$.
	Since $F-H=2iU$ on $\bbC$
	and $U(z)=U'(z)=0$, we have
	\[
	F(z)=H(z),\qquad
	F'(z)=H'(z).
	\]
	In this case, take $\tau_z=1$.
	
	Finally, by the unitarity of the Bargmann transform,
	there exist $f,h\in L^2(\bbR)$ such that
	\[
	\mathcal Bf=F,\qquad
	\mathcal Bh=H.
	\]
	Since $F\not\sim H$, we have $f\not\sim h$.
	For every $z\in\Lambda$,
	we have $\overline z\in\Lambda$, and hence
	\[
	F(\overline z)=\tau_{\overline z}H(\overline z),
	\qquad
	F'(\overline z)=\tau_{\overline z}H'(\overline z).
	\]
	By \eqref{eq:hermite},
	\[
	V_{h_0}f(z)=\tau_{\overline z}V_{h_0}h(z),
	\qquad
	V_{h_1}f(z)=\tau_{\overline z}V_{h_1}h(z).
	\]
	Let $g=a h_0+b h_1$, where $a,b\in\bbC$. Then
	\begin{align*}
		V_gf(z)
		&=
		\overline a\,V_{h_0}f(z)
		+\overline b\,V_{h_1}f(z)\\
		&=
		\tau_{\overline z}
		\left(
		\overline a\,V_{h_0}h(z)
		+\overline b\,V_{h_1}h(z)
		\right)\\
		&=
		\tau_{\overline z}V_gh(z).
	\end{align*}
	Since $|\tau_{\overline z}|=1$,
	\[
	|V_g f(z)| = |V_g h(z)|,\qquad
	\forall g\in\mathcal H_{01},
	\forall z\in\Lambda.
	\]
	
	Now let $v=(x_0,\omega_0)\in\bbR^2$ and define
	\[
	f_v(t)=\ee{2\pi i\omega_0t}f(t-x_0),\qquad
	h_v(t)=\ee{2\pi i\omega_0t}h(t-x_0).
	\]
	Translation and modulation are unitary, so $f_v,h_v\in L^2(\bbR)$ and $f_v\not\sim h_v$.
	A change of variables in \eqref{eq:stft} gives
	\[
	V_gf_v(x,\omega)
	=\ee{-2\pi i(\omega-\omega_0)x_0}
	V_gf(x-x_0,\omega-\omega_0),
	\]
	and the same identity holds for $h_v$. Then we have
	\[
	|V_g f_v(z)| = |V_g h_v(z)|,\qquad
	\forall g\in\mathcal H_{01},
	\forall z\in v+\Lambda.
	\]
	This proves the sharpness assertion for every translated rectangular lattice in Theorem~\ref{thm:lattice}.
\end{proof}

For example, the four windows in \eqref{eq:fourwindows} do phase retrieval on $2^{-1/2}\bbZ\times2^{-1/2}\bbZ$.
If the spacing in both directions is increased to $\ell>2^{-1/2}$, then not only these four windows but also all windows in $\mathcal H_{01}$ fail to distinguish the two signals constructed in the proof of Theorem~\ref{thm:lattice} on $\ell\bbZ^2$.

\textbf{Conflict of interest statement}.
The authors state that there is no conflict of interest.

\textbf{Data availability statement}.
The authors state that the paper contains no data.


\end{document}